\documentclass{amsart}
\usepackage{latexsym}
\usepackage{graphics}
\usepackage{amssymb}
\usepackage[mathscr]{eucal}
\usepackage{xspace}
\usepackage{graphicx}

\newcommand{\af}{\alpha}
\newcommand{\bt}{\beta}
\renewcommand{\d}{\delta}
\newcommand{\D}{\Delta}
\newcommand{\e}{\varepsilon}

\newcommand{\ep}{\varepsilon}
\newcommand{\gm}{\gamma}
\newcommand{\Gm}{\Gamma}
\renewcommand{\l}{\lambda}

\renewcommand{\o}{\omega}
\renewcommand{\O}{\Omega}
\renewcommand{\S}{\sigma}

\newcommand{\R}{{\mathbb R}}
\newtheorem{theorem}{Theorem}[section]
\newtheorem{lemma}{Lemma}[section]
\newtheorem{proposition}[lemma]{Proposition}
\newtheorem{corollary}{Corollary}[section]
\newtheorem{definition}{Definition}[section]
\newtheorem{remark}{Remark}[section]
\begin{document}
\baselineskip 15pt

\title[Well-posedness for steady supersonic flows past a Lipschitz wedge]
{Well-posedness for two-dimensional steady supersonic Euler flows
past a Lipschitz wedge }
\author{Gui-Qiang Chen $\qquad$  Tian-Hong Li}
\address{G.-Q. Chen, Department of Mathematics \\Northwestern University
\\Evanston, Illinois 60208, USA}
\email{gqchen@math.northwestern.edu}
\address{T.-H Li, Academy
of Mathematics and Systems Science, Chinese Academy of Sciences,
Beijing 100080, PRC} \email{thli@math.ac.cn} 
.
\medskip
\date{September 1, 2006}

\begin{abstract}
For a supersonic Euler flow past a straight wedge whose vertex angle
is less than the extreme angle, there exists a shock-front emanating
from the wedge vertex, and the shock-front is usually strong
especially when the vertex angle of the wedge is large. In this
paper, we establish the $L^1$ well-posedness for two-dimensional
steady supersonic Euler flows past a Lipschitz wedge whose boundary
slope function has small total variation,
when the total variation of the incoming flow is sufficiently small.
In this case, the Lipschitz wedge perturbs the flow and the waves
reflect after interacting with the strong shock-front or the wedge
boundary. We first obtain the existence of solutions in $BV$ when
the incoming flow has small total variation by the wave front
tracking method and then study the $L^1$ stability of the solutions.
In particular, we incorporate the nonlinear waves generated from the
wedge boundary to develop a Lyapunov functional between two
solutions, which is equivalent to the $L^1$ norm, and
prove that the functional decreases in the flow direction. Then the
$L^1$ stability is established, so is the uniqueness of the
solutions by the wave front tracking method. Finally, we show the
uniqueness of solutions in a broader class, i.e. the class of
viscosity solutions.
\end{abstract}
\maketitle

\numberwithin{equation}{section}
\setcounter{page}{1}

\section {Introduction}
\setcounter{equation}{0}
For the Cauchy problem of a strictly
hyperbolic system of conservation laws:
\begin{equation}\label{1.1}
U_t + F(U)_x=0, \qquad U\in \R^n,
\end{equation}
\begin{equation}\label{1.2}
U|_{t=0}=\overline{U}(x),
\end{equation}
whose each characteristic field is either linearly degenerate or
genuinely nonlinear, the existence of weak solutions to
(\ref{1.1})--(\ref{1.2}) with small total variation was first proved
by Glimm \cite{G} by a probabilistic algorithm, the Glimm scheme;
and a deterministic version of the Glimm scheme was developed by Liu
\cite{Liu}. Alternative methods for constructing solutions of the
Cauchy problem were first introduced in \cite{Da,Di}, based on wave
front tracking. For the scaler equation, $F$ is approximated by
piecewise linear functions $F_{\nu}$ in Dafermos \cite{Da} so that
the approximate solutions are piecewise constants and all the
interactions are determined by solving the Riemann problem. The
method was generalized to the $2\times 2$ case in DiPerna \cite{Di}
in which piecewise constant approximate solutions are constructed so
that the wave interactions can be determined by only solving the
Riemann problem. Bressan \cite{Br} developed the wave front tracking
method for $n\times n$ systems by overcoming the difficulty that the
procedure used in \cite {Di} may yield an infinite number of
discontinuities in finite time when $n>2$; and the wave front
tracking method was further simplified later in Baiti-Jenssen
\cite{BJ}. Also see Bressan \cite{Bressan-book}, Dafermos
\cite{Dafermos-book}, Holden-Risebro \cite{HR-book}, and LeFloch
\cite{lefloch-book} for further references.

Within the class of initial data $\overline{U}\in L^1\cap BV(R;
R^n)$ with suitably small total variation, it was established that
problem (\ref{1.1})--(\ref{1.2}) is well-posed in $L^1(R; R^n)$ for
the solutions generated by the wave front tracking algorithm. In
particular, in Bressan-Colombo \cite{BC}, Bressan-Crasta-Piccoli
\cite{BCP}, and Bressan-Liu-Yang \cite{BLY}, it was proved that the
entropy solutions of (\ref{1.1})--(\ref{1.2}) constitute a semigroup
which is Lipschitz continuous with respect to time and initial data.
Lewicka-Trivisa \cite{LT} obtained the $L^1$ well-posedness of
solutions generated by the wave front tracking method for the Cauchy
problem (\ref{1.1})--(\ref{1.2}) with the initial data
$\overline{U}$ being a small perturbation of a fixed Riemann problem
$(U_-, U_+)$ containing two large shocks, under the necessary
stability condition (cf. \cite{BC2,LT,Le2}; also see
\cite{Le1,Le3}). The $L^1$ well-posedness in the class of viscosity
solutions for the Cauchy problem has been also established (cf.
\cite{Br2,BiB,C} and the references therein).


In this paper, we are concerned with the $L^1$ well-posedness of a
physical nonlinear problem of initial-boundary value type, which
governs two-dimensional steady supersonic Euler flows past a curved
wedge.
More specifically, the two-dimensional steady supersonic Euler flows
are generally governed by
\begin{equation}\label{eq1.1}
\left\{%
\begin{array}{ll}
    (\rho u)_x + (\rho v)_y = 0 , \\
    (\rho u^2 + p)_x + (\rho uv)_y = 0,  \\
   (\rho uv)_x + (\rho v^2 + p)_y = 0,\\
   (\rho u(E+p/\rho))_x+ (\rho v(E+p/\rho))_y=0,
  \end{array}%
\right.
\end{equation}
where $(u,v)$ is the velocity, $\rho$ the density, $p$ the scalar
pressure, and $ E=\frac{1}{2}(u^2+v^2)+e(\rho,p)$ the total energy
with $e$ the internal energy (a given function of $(\rho,p)$ defined
through thermodynamical relationships). The other two thermodynamic
variables are the temperature $\theta$ and the entropy $S$.
For an ideal gas,
\begin{equation}
p=R\rho \theta, \qquad e=c_v\theta, \qquad \gamma=1+\frac{R}{c_v}>1,
\label{1.5a}
\end{equation}
and
\begin{equation}
p=p(\rho,S)=\kappa\rho^\gamma e^{S/c_v}, \qquad
e=\frac{\kappa}{\gamma-1}\rho^{\gamma-1}e^{S/c_v}
   =\frac{R\theta}{\gamma-1},
\label{1.6a}
\end{equation}
where $R, \kappa,$ and $c_v$ are all positive constants.

If the flow is isentropic, i.e. $S=const.$, then the pressure $p$ is
a function of the density $\rho$, $p=p(\rho)$, and the flow is
governed by the following simpler isentropic Euler equations:
\begin{equation}\label{eq1.2}
\left\{%
\begin{array}{ll}
    (\rho u)_x + (\rho v)_y = 0 , \\
    (\rho u^2 + p)_x + (\rho uv)_y = 0,  \\
   (\rho uv)_x + (\rho v^2 + p)_y = 0.\\
  \end{array}%
\right.
\end{equation}
For polytropic isentropic gases, by scaling,
the pressure-density relationship can be expressed as
\begin{equation}\label{1.2b}
p(\rho)=\rho^\gamma/\gamma, \qquad \gamma>1.
\end{equation}
For the isothermal flow, $\gamma=1$.
The quantity
$$
c=\sqrt{p_\rho(\rho, S)}
$$
is defined as the sonic speed and, for polytropic gases,
$c=\sqrt{\gamma p/\rho}$.

System (\ref{eq1.1}) or (\ref{eq1.2}) governing a supersonic flow
(i.e., $u^2+v^2>c^2$) has all real eigenvalues and is hyperbolic,
while system (\ref{eq1.1}) or (\ref{eq1.2}) governing a subsonic
flow (i.e., $u^2+v^2<c^2$) has complex eigenvalues and
is elliptic-hyperbolic mixed and composite.

\medskip
The study of two-dimensional steady supersonic flows past a wedge
can date back to the 1940s (cf. Courant-Friedrichs \cite{CF}). Local
solutions around the wedge vertex were first constructed in Gu
\cite{Gu}, Li \cite{Li}, Schaeffer \cite{Schaeffer}, and the
references cited therein. Global potential solutions were
constructed in various different setups in
\cite{Chen1,Chen2,Chen3,Zh2}
when the wedge vertex angle is less than the critical angle.

For the full Euler equations, when a wedge is straight and the wedge
vertex angle is less than the critical angle, there exists a
supersonic shock-front emanating from the wedge vertex so that the
constant states on both sides of the shock-front are supersonic; the
critical angle condition is necessary and sufficient for the
existence of the supersonic shock
(cf. Courant-Friedrichs \cite{CF}). When the incoming flow is
uniform,  Chen-Zhang-Zhu \cite{CZZ} first established the existence
of global supersonic Euler flows, especially the nonlinear
structural stability of the strong shock-front emanating from the
wedge vertex under the $BV$ perturbation of the wedge boundary.
In this paper, we first show the existence of solutions to the above
problem when the incoming flow is a $BV$ perturbation of the uniform
flow by the wave front tracking method, and then we establish the
$L^1$-stability of entropy solutions generated by this method. Based
on these, we establish estimates on the uniformly Lipschitz
semigroup $\mathcal{P}$ generated by the wave front tracking limit
and prove the uniqueness of solutions by means of local integral
estimates within a broader class of solutions, i.e. the class of
viscosity solutions.

One of the main new ingredients in our approach here is to develop
techniques to handle the boundary difficulty, in comparison with the
earlier works on the Cauchy problem. For the $L^1$ stability of
solutions of the Cauchy problem, the decrease of the Lyapunov
functional was achieved by essentially using the cancellation of the
distances on both sides of waves. However, for our wedge problem
that is the problem of initial-boundary value type,
there is no such cancellation near the boundary, since only one-side
is possible near the wedge boundary. In order to overcome this
difficulty, we employ the exact feature of the boundary condition to
obtain an estimate and refine the functional based on this estimate.
In particular, since the flow of two solutions near the boundary
must be parallel, we identify the relation between these two states,
which is desirable for redesigning the functional to ensure the
decreasing of the functional in the flow direction.

For concreteness, as in Chen-Zhang-Zhu \cite{CZZ}, we will analyze
the problem in the region below the lower side $\Gamma$ of the wedge
for the Euler flows for $U=(u,v,p,\rho)$ governed by system
(\ref{eq1.1}) and $U=(u,v,\rho)$ by (\ref{eq1.2}); the case above
the wedge can be handled in the same fashion. Then we have

\begin{enumerate}\renewcommand{\theenumi}{\roman{enumi}}
\item There exists a Lipschitz function $g\in Lip(\R_+)$
with $g'\in BV(\R_+)$, $g'(0+)=0$, and $g(0)=0$ such that
$$
\quad \Omega:=\{ (x,y)\,\, : \,\, y< g(x),\, x\ge 0 \}, \qquad
 \Gamma:= \{ (x,y)\,\,:\,\, y= g(x),\, x \ge 0 \},
$$
and
${\bf n}(x\pm) =\frac{(-g'(x\pm),1)}{\sqrt{(g'(x\pm))^2+1}}$
is the outer normal vector to $\Gamma$ at the point $x\pm$
(see Fig. \ref{fig:1});

\begin{figure}[htbp] \begin{center}
\includegraphics[height=1.5in,width=3in]{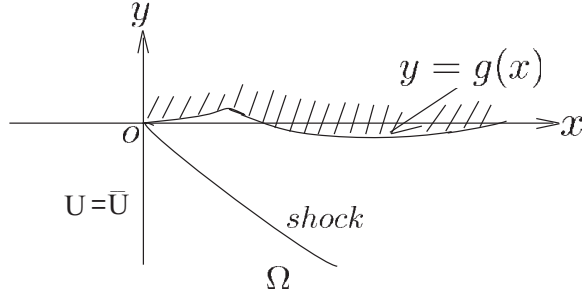}
\caption{Supersonic flow past a curved wedge} \label{fig:1}
\end{center} \end{figure}

\item The upstream flow $U=\bar{U}(y)=(\bar{u}, \bar{v}, \bar{p}, \bar{\rho})(y)$
at $x=0$ satisfies
$$
\bar{u}(y)>0, \quad \bar{v}(y)>0, \quad \bar{u}(y)^2+\bar{v}(y)^2 >
\bar{c}(y)^2:=\gamma \bar{p}(y)/\bar{\rho}(y),
$$
and
$$
0< \arctan(\bar{v}(y)/\bar{u}(y)) < \omega_{crit},
$$
where $\omega_{crit}$ is the critical vertex angle so that there is
a supersonic shock-front emanating from the wedge vertex.
\end{enumerate}

With this setup, the wedge problem can be formulated into the
following problem of initial-boundary value type for system
(\ref{eq1.1}) or (\ref{eq1.2}):

{\it Cauchy Condition:}
\begin{equation}\label{IC}
U|_{x=0} = \overline{U}(y);
\end{equation}

{\it Boundary Condition}:
\begin{equation}\label{BC}
(u, v) \cdot {\bf n}=0\qquad\mbox{ on }\, \Gamma.
\end{equation}

This paper is organized as follows. In Section 2, we discuss the
basic properties for the adiabatic Euler equations and related
nonlinear waves. In Section 3, we discuss the wave front tracking
method and define the interaction potential $Q$, and then we prove
the existence of entropy solutions to the initial-boundary value
problem. In Section 4, we construct the Lyapunov functional $\Phi$
between two solutions to incorporate the nonlinear waves generated
by the wedge boundary vertices, which is equivalent to the $L^1$
distance between these two solutions. In Section 5, we prove the
decrease of $\Phi$ in the flow direction, which implies the $L^1$
stability of the solutions. In Section 6, we prove the existence of
the semigroup generated by the wave front tracking method and
establish some estimates on the uniformly Lipschitz semigroup
$\mathcal{S}$ generated by the wave front tracking limit. In Section
7, we prove the uniqueness of solutions by means of local integral
estimates within a broader class of solutions.

\section{Euler Equations and Nonlinear Waves}
\setcounter{equation}{0} In this section, we review some basic
properties of the adiabatic Euler equations (\ref{eq1.1}) and
related nonlinear waves, which will be used in the subsequent
sections. The Euler equations for steady supersonic flows can be
written in the following conservation form:
\begin{equation}\label{6.1.1}
 W(U)_x + H(U)_y = 0,
\qquad\,\, U=(u, v, p,\rho)^\top,
\end{equation}
with
$$
W(U)=(\rho u,\rho u^2+p,\rho uv, \rho u(h+\frac{u^2+v^2}{2}))^\top,
 H(U)=(\rho v,\rho uv,\rho v^2+p, \rho
v(h+\frac{u^2+v^2}{2}))^\top,
$$
and  $h=\frac{\gamma p}{(\gamma-1)\rho}$. The eigenvalues of system
(\ref{6.1.1}) are
\begin{equation}\label{eq6.1.3}
\l_j=\frac{uv+(-1)^jc\sqrt{u^2+v^2-c^2}}{u^2-c^2},\;\; j=1,4; \qquad
\l_{i}=v/u, \,\, i=2,3,
\end{equation}
where $c^2 = \gm p/\rho$. If the flow is supersonic (i.e.
$u^2+v^2>c^2$), system \eqref{6.1.1} is hyperbolic; and, in
particular, when $u>c$, system \eqref{6.1.1} has the four
corresponding linearly independent eigenvectors:
\begin{eqnarray}\label{2.5a}
&& {\bf r}_j={\bf\kappa}_j
 (-\l_j, 1, \rho (\l_j u-v), \rho (\l_j u-v)/c^2)^\top, \quad j=1,4,\\
&&\nonumber {\bf r}_2= (u, v, 0, 0)^\top, \qquad {\bf
r}_3=(0,0,0,\rho)^\top,
\end{eqnarray}
where ${\bf\kappa}_j$ are chosen so that ${\bf
r}_j\cdot\nabla\l_j=1$ since the $j$th-characteristic fields are
genuinely nonlinear, $j=1,4$. Note that the second and third
characteristic fields are always linearly degenerate: ${\bf
r}_j\cdot \nabla \lambda_j=0, j=2,3. $

\begin{definition}[Entropy Solutions]\label{D.6.1}
A $BV$ function $U=U(x,y)$ is called an entropy solution of the
initial-boundary value problem {\rm (\ref{6.1.1})} and {\rm
(\ref{IC})}--{\rm (\ref{BC})} provided that

{\rm (i)} $U$ is a weak solution of {\rm (\ref{6.1.1})} and
satisfies
$$
U|_{x=0}=\bar{U}(y), \quad (u,v)\cdot {\bf n}|_{y=g(x)}=0
\qquad\mbox{in the trace sense};
$$

{\rm (ii)} $U$ satisfies the entropy inequality, i.e. the steady
Clausius inequality:
\begin{equation}\label{entropy-ineq:6.1}
(\rho u S)_x +(\rho v S)_y\ge 0
\end{equation}
in the sense of distributions in $\Omega$ including the wedge
boundary.
\end{definition}

We now discuss the wave curves in the phase space. The contact
Hugoniot curves $C_i(U_0)$ through $U_0$ are
\begin{equation}\label{6.3.5}
C_i(U_0):\quad p= p_0,\quad w = v/u = v_0/u_0,\quad\, i=2,3,
\end{equation}
which describe compressible vortex sheets. We remark that, although
the two contact discontinuities coincide as a single vortex sheet in
the physical $xy$-plane,  it requires two independent parameters to
describe them in the phase space $U=(u,v,p,\rho)$ since there are
two linearly independent eigenvectors corresponding to the repeated
eigenvalues $\lambda_2=\lambda_3=v/u$ of the two linearly degenerate
fields.

Moreover, the rarefaction wave curves $R_j(U_0)$
in the phase space through $U_0$ are
\begin{equation}\label{6.3.6}
R_j^-(U_0):\,\, dp = c^2d\rho,\,\, du = - \lambda_j dv,\,\,
\rho(\lambda_j u - v)dv = dp \qquad\mbox{for}\,\, \rho <\rho_0,
\qquad \hbox{$j = 1,4$.}
\end{equation}

The Rankine-Hugoniot conditions for (\ref{6.1.1}) are
\begin{eqnarray}
s[W(u)] = [H(u)],
     \label{eq6.1.8a}
\end{eqnarray}
where $s$ is the propagation speed of the discontinuity. Then
$$
s=s_j:=\frac{u_0v_0 +(-1)^j \overline{c}_0\sqrt{u_0^2+v_0^2-
\overline{c}_0^2}} {u_0^2-\overline{c}_0^2}, j=1,4, \quad \S =
\S_{i} = v_0/u_0, i=2,3,
$$
where $\overline{c}_0^2 = \frac{c_0^2}{b_0} \frac{\rho}{\rho_0}$ and
$b_0= \frac{\gm +1}{2}-\frac{\gm -1}{2} \frac{\rho}{\rho_0}$.

Plugging $s_i, i=2,3,$ into (\ref{eq6.1.8a}), we obtain the
$i^{th}$-contact Hugoniot curves $C_i(U_0), i=2,3, $ as defined in
(\ref{6.3.5}); while plugging $s_j, j=1,4,$ into (\ref{eq6.1.8a}),
we obtain the $j^{th}$-Hugoniot curve $S_j(U_0)$ through $U_0$:
\begin{equation*}
S_j(U_0):\,\, [p]= \frac{c_0^2}{b_0}[\rho], \,\, [u]= -s_j [v],\,\,
\rho_0(s_j u_0 - v_0)[v]=[p] \qquad\mbox{for}\,\, \rho>\rho_0,
\qquad \hbox{j = 1, 4}.
\end{equation*}
The half curves of $S_j(U_0)$ for $\rho>\rho_0$, denoted by
$S_j^+(U_0), j=1,4$, in the phase space are called the shock curves
on which any state with $U_0$ forms a shock in the $x-y$ plane
satisfying the entropy condition as explain in Lemma 2.1 below.

Note that $S_j^+(U_0)$ contacts with $R_j^-(U_0)$ at $U_0$ up to
second-order.

As indicated in \cite{CZZ}, we have

\begin{lemma} If $U$ is a piecewise smooth solution, then, on the shock
wave, the entropy inequality {\rm (\ref{entropy-ineq:6.1})} in
Definition {\rm 2.1} is equivalent to any of the following:

{\rm (i)} The physical entropy condition: the density increases
across the shock in the flow direction,
\begin{equation}\label{entropy-ineq:2a}
\rho_{front}<\rho_{back};
\end{equation}

{\rm (ii)} The Lax entropy condition: on the $j^{th}$-shock with the
shock speed $\sigma_j$,
\begin{eqnarray}
&&\lambda_j(back)<s_j<\lambda_j(front), \qquad j=1,4,\\
&& s_1<\lambda_{2,3}(back), \qquad \lambda_{2,3}(front)<s_4.
\end{eqnarray}
\end{lemma}

The following properties and related estimates of wave interactions
in Lemmas 2.2--2.8 have been obtained in Chen-Zhang-Zhu \cite{CZZ}.
We list them below for subsequent use in this paper.

\subsection{Riemann problems and Riemann solutions}
\label{sec:19} We start with Riemann problems and their solutions.

{\it Lateral Riemann problem.}
The simplest case of problem (\ref{6.1.1}) and
(\ref{IC})--(\ref{BC}) is $g\equiv 0$. It can be shown that, if
$g\equiv 0$, then problem (\ref{6.1.1}) admits an entropy solution
that consists of a constant state $U_-$ and a constant state $U_+$,
with $U_+=(u_+,0,p_+,\rho_+)$ and $u_+ >c_+>0$ in the subdomain of
$\Omega$ separated by a straight shock emanating from the vertex.
That is to say that the state ahead of the shock-front is $U_-$,
while the state behind the shock-front is $U_+$.
When the angle between the flow direction of the front state and the
wedge boundary at a boundary vertex is larger than $\pi$, then an
entropy solution contains a rarefaction wave that separates the
front state from the back state.

{\it Riemann problem involving only weak waves.}
Consider the following initial value problem:
\begin{equation} \label{eq6.2.1}
\left\{%
\begin{array}{ll}
    W(U)_x + H(U)_y = 0,  \\
    U|_{x=x_0}= \underline{U}= \left\{%
\begin{array}{ll}
    U_a, & \hbox{$y > y_0$,} \\
    U_b, & \hbox{$y < y_0$,}
\end{array}%
\right.
\end{array}%
\right.
\end{equation}
where
$U_b$ and $U_a$ are constant states.

\begin{lemma}\label{L:6.2.1}
There exists $\e>0$ such that, for any states $U_a,U_b\in
O_{\e}(U_+)$ or $U_a,U_b\in
O_{\e}(U_-)$, problem {\rm (\ref{eq6.2.1})} admits a unique
admissible solution consisting of four elementary waves.
\end{lemma}

{\it Riemann problem involving a strong 1-shock.}
For simplicity, we use notation $\{U_b,U_a \} =(\af_1,\af_2,\af_3,
\af_4)$ to denote the solution to the Riemann problem, where $\af_i$
is the strength of  the $ith$ elementary wave. For any $U\in
S_1(U_-)$, we also use $\{ U_-, U\}=(\S,0,0,0)$ to denote the
$1$-shock that connects $U_-$ and $U$ with speed $\S$. Then we have

\begin{lemma}\label{L:6.2.3}
Let $\{U_-,U_+\}=(\S_0,0,0,0),\rho_+>\rho_-,$ and $\gm >1$.
Then
$$
\S_0 < 0, \qquad u_+ <u_- <(1+1/\gamma)u_+,
$$
and
$$
{\rm det}(\nabla_U H(U_+)-\S_0\nabla_U W(U_+))>0.
$$
Furthermore,
there exists a neighborhood $O_\e(U_+)$ of $U_+$ and a neighborhood
$O_\e(U_-)$ of $U_-$ such that $U_0\in O_\e(U_-)$ and the shock
polar $S_1(U_0)\cap O_\e(U_+)$ can be parameterized by the shock
speed $\sigma$ as $ \S \rightarrow G(U_0, \S)
$
with $G \in C^2$ near
$(U_-, \S_0)$ and $G(U_-, \S_0)= U_+$.
\end{lemma}

\subsection{Estimates on wave interactions and reflections}
\label{sec:24}
We have
\begin{lemma}[Estimates on weak wave
interactions]\label{prop6.3.1} Suppose that $U_b, U_m, U_a\in
O_\e(U_+)$, or  $U_b, U_m, U_a\in O_\e(U_-)$, are three states with
$\{U_b, U_m\}=(\af_1, \af_2, \af_3, \af_4)$, $\{U_m, U_a\} =(\bt_1,
\bt_2, \bt_3, \bt_4)$, and $\{U_b, U_a\}=(\gm_1, \gm_2, \gm_3,
\gm_4) $.
Then
$$
\gm_i=\af_i+\bt_i+O(1) \triangle(\af, \bt),
$$
where
$
\triangle(\af, \bt)=(|\af_4| +|\af_3|+|\af_2|)|\bt_1|
  + |\af_4|(|\bt_2| + |\bt_3|) +
\sum_{j=1,4}\triangle_j(\af, \bt)
$
with
$$
\triangle_j(\af, \bt)=\left\{
\begin{array}{ll}
    0, & \hbox{$\af_j \ge 0, \; \bt_j \ge 0$,} \\
    |\af_j||\bt_j|, & \hbox{otherwise.}
\end{array}%
\right.
$$
\end{lemma}


\medskip Denote $\{C_k(a_k, b_k)\}_{k=0}^\infty$ the points
$\{(a_k,b_k)\}_{k=0}^\infty$ in the $xy$-plane with $a_{k+1}>a_k>0$.
Set
\begin{eqnarray}
&&\omega_{k, k+1}=\arctan\left(\frac{b_{k+1}-b_k}{a_{k+1}-a_k}\right),
\quad
\omega_k =\omega_{k, k+1} - \omega_{k-1,k}, \quad \omega_{-1,0}=0,
\nonumber\\
&&\Omega_{k+1}=\{(x,y)\, :\,  x \in [a_k, a_{k+1}),\;
y<b_k+(x-a_k)\tan(\omega_{k,k+1})\},
\label{3.ref-1}\\
&&\Gamma_{k+1}=\{(x,y)\, :\, x \in [a_k, a_{k+1}),\;
 y = b_k+(x-a_k)\tan(\omega_{k,k+1})\},\nonumber
\end{eqnarray}
and the outer normal vector to $\Gamma_k$:
\begin{equation}\label{3.ref-2}
{\bf n}_{k+1}
=\frac{(b_{k}-b_{k+1},a_{k+1}-a_k)}{\sqrt{(b_{k+1}-b_k)^2+(a_{k+1}-a_k)^2}}
=(-\sin(\omega_{k,k+1}), \cos(\omega_{k, k+1})).
\end{equation}


Then we consider the initial-boundary value problem
with $\underline{U}$ a constant state:
$$
\left\{%
\begin{array}{ll}
    (\ref{6.1.1}) & \hbox{ in $\O_{k+1}$}, \\
    U|_{x=a_{k}}=\underline{U}, & \hbox{} \\
    (u,v)\cdot {\bf n}_{k+1}=0 & \hbox{on $\Gm_{k+1}$}.
\end{array}%
\right.
$$

\begin{lemma}[Estimate on the boundary perturbation of
the strong shock]\label{prop6.2} For $\e>0$ sufficiently small,
there exists $\hat{\e}=\hat{\e}(\e)<\e$ so that
$G(O_{\hat{\e}}(\S_0))\subset O_\e(U_+)$ and, when $|\o_k|<\e$, the
equation
$G(\S) \cdot ({\bf n}_k, 0, 0) =0$
admits a unique solution $\S_k \in O_{\hat{\e}}(\S_0)$.
Moreover, we have
\begin{equation}\label{eq6.3.21}
\S_{k+1}=\S_k+ K_{bs} \o_k,
\end{equation}
where $|K_{bs}|$ is bounded.
\end{lemma}
\begin{lemma}[Estimate on the boundary perturbation of
weak waves] Let $U_k=(u_k,v_k,p_k,\rho_k)$ be the  state near the
boundary with $ (u_k,v_k)\cdot {\bf n}_k=0. $ Then there exists
$U_{k+1}$ such that $ \{ U_k, U_{k+1}\}=(\delta_1, 0, 0, 0)$ and
$(u_{k+1}, v_{k+1})\cdot{\bf n}_{k+1}=0$. Furthermore,
$$
\delta_1=K_{b0}\o_k,
$$ where $K_{b0}$ is bounded.
\end{lemma}

\begin{lemma}[Estimates on the reflection  of weak waves
on the boundary] Let $\{U_b, U_k\}=(0, \af_2, \af_3, \af_4)$ and
$(u_k,v_k)\cdot {\bf n}_k=0$.  Then there exists $U_{k+1}$ such that
$ \{ U_b, U_{k+1}\}=(\delta_1, 0, 0, 0)$ and $(u_{k+1},
v_{k+1})\cdot{\bf n}_{k}=0$. Furthermore,
$$
\delta_1=K_{b4}\af_4+K_{b3}\af_3+K_{b2}\af_2,
$$
where $K_{b4}, K_{b3}, K_{b2}$, and $K_{b0}$ are $C^2-$functions of
$(\af_4, \af_3, \af_2, \bt_1, \o_{k}; U_b)$ satisfying
\begin{eqnarray*}
&&K_{b4}|_{\{\o_{k}=\af_4=\af_3=\af_2=\bt_1=0, U_b=U_+\}}=1,\\
&&K_{b2}|_{\{\o_{k}=\af_4=\af_3=\af_2=\bt_1=0, U_b=U_+\}}=
K_{b3}|_{\{\o_{k}=\af_4=\af_3=\af_2=\bt_1=0, U_b=U_+\}}=0.
\end{eqnarray*}
\end{lemma}


\begin{lemma}[Estimates on the interaction between the strong shock
and weak waves from above]\label{prop6.3.3} Let $U_m,U_a\in
O_\e(U_+)$ with $ \{G(U_b, \S),U_m\}=(0, 0, 0, 0)$ and
$\{U_m,U_a\}=(\bt_1, 0, 0, 0)$.  Then there exists a unique
$(\S',\d_2,\d_3,\d_4)$ such that the Riemann problem {\rm
(\ref{eq6.2.1})} with $U_b\in O_\e(U_-)$ admits an admissible
solution consisting of a strong $1$-shock, two contact
discontinuities of strengths $\d_2$ and $\d_3$, and a weak $4$-wave
of strength $\d_4$:
$$
\{U_b,U_a\}=(\S',\d_2, \d_3, \d_4).
$$
Moreover,
$$
\S'=\S+ K_{s1}\bt_1, \quad\; \d_2=K_{s2}\bt_1, \quad\,
\d_3=K_{s3}\bt_1, \quad\; \d_4=K_{s4}\bt_1,
$$
where $ |K_{s4}|<1$, and $|K_{s1}|+|K_{s2}|+|K_{s3}|$ is bounded.
Furthermore, we have
\begin{equation}\label{key}
\widetilde{K_{s4}}|\frac{\l_{4+}-\S_0}{\l_{1+}-\S_0}| = |\frac{\S_0
u_-Q - u_+ \l_{4+} P}{\S_0 u_-Q + u_+\l_{4+}P}|<1.
\end{equation}
\end{lemma}


\begin{lemma}[Estimates on the interaction between
the strong shock and weak waves from below]\label{prop6.3.4} Let
$U_m,U_b\in O_\e(U_-)$ and $U_a\in O_\e(U_+)$ with
$$
\{U_b,U_m\}=(\af_1, \af_2, \af_3, \af_4),\qquad
\{U_m, U_a\}=(\S, \bt_2, \bt_3, \bt_4).
$$
Then there exists a unique $(\S',\d_2,\d_3,\d_4)$ such that the
Riemann problem {\rm (\ref{eq6.2.1})}  admits an
admissible solution consisting of a strong $1$-shock, two contact
discontinuities of strengths $\d_2$ and $\d_3$, and a weak $4$-wave
of strength $\d_4$:
$$
\{U_b,U_a\}=(\S',\d_2, \d_3, \d_4).
$$
Moreover,
$$
\S'=\S+ \sum_{i=1}^4K_{1i}\af_i+O(1)\D, \quad\;
\d_2= \bt_2+\sum_{i=1}^4 K_{2i}\af_i+O(1)\D,
$$
$$
\d_3=\bt_3+\sum_{i=1}^4 K_{3i}\af_i+O(1)\D,
\quad\;\d_4=
\bt_4+\sum_{i=1}^4 K_{4i}\af_i+O(1)\D,
$$
where $|K_{ji}|, i, j=1,..,4,$ are bounded and
$\Delta=\sum_{i=1,2,3,4, j=2,3,4}|\af_i\bt_j|.$ Furthermore, we can
write the estimates in a more precise fashion:
\begin{eqnarray*}
\sigma'=\sigma+\sum_{i=1}^4 \widetilde{K_{1i}}\af_i, \,\,
\delta_2=\bt_2+\sum_{i=1}^4  \widetilde{K_{2i}}\af_i,\,\,
\delta_3=\bt_3+ \sum_{i=1}^4 \widetilde{K_{3i}}\af_i,\,\,
\delta_4=\bt_4+ \sum_{i=1}^4 \widetilde{K_{4i}}\af_i,
\end{eqnarray*}
where $ \sum_{i,j=1}^4|\widetilde{K_{ji}}|\le M$ for some $M>0$.
\end{lemma}
\begin{proof}
We first consider the interaction between $(U_b,
U_m)=(\af_1,\af_2,\af_3,\af_4)$ and $(U_m, G(U_m, \S))=(\S, 0,0,0)$
to find that the solution is the perturbation of the unperturbed
states of the strong shock.  From Lemma \ref {L:6.2.3}, we know that
$U=G(U_0, \S) $ near $(U_-, \S_0)$ with $G\in C^2$, which implies
that $U =(u, v, p, \rho)$ depends continuously on the state
$U_0=(u_0, v_0, p_0, \rho_0)$.  The perturbation near $U_-$ is
equivalent to $\sum_{i=1}^4O(1)\af_i$. Then the interaction estimate
between  $(U_b, G(U_m, \S))$ and $(G(U_m, \S), U_a)$ follows from
Lemma \ref{prop6.3.1}.
\end{proof}

\section{Wave front tracking method and existence of entropy solutions}
\setcounter{equation}{0}


The basic idea of the wave front tracking method is to construct
approximate solutions within a class of piecewise constant
functions: First, approximate the initial data by a piecewise
constant function and solve the resulting Riemann problems and
replace the rarefaction waves by step functions with many small
discontinuities; then track the outgoing fronts until the first time
when two waves interact which are determined by a new Riemann
problem; and finally design a simplified Riemann solver so that the
number of wave fronts is finite for all $x\ge 0$ in the flow
direction.

\subsection{The Riemann solvers}
\setcounter{equation}{0} As mentioned in Section 2, the solution to
the Riemann problem $(U_b, U_a)$ is  a self-similar solution given
by at most five states separated by shocks or rarefaction waves.
More precisely, there exists $C^2$ curves
$\alpha\rightarrow\psi(\alpha)(u)$ parameterized by arc length such
that
$$
U_b=\psi_4(\alpha_4)\circ\ldots\circ\psi_1(\alpha_1)(U_a)
$$
for some $\alpha=(\alpha_1, \dots, \alpha_4)$ and
$U_i=\psi_i(\alpha_i)\circ\ldots\circ\psi_1(\alpha_1)(U_a)$. When
$\alpha_i$ is positive (negative), states $U_{i-1}$ and $U_i$ are
separated by an $i$-rarefaction ($i$-shock) wave, so we call
$\alpha_i$ the strength of the $i$-wave.

For given initial data $\overline{U}$, let $\overline{U}^\e, \e>0,$
be a sequence of piecewise constant functions approximating
$\overline{U}$ in the $L^1$ norm, and the wedge boundary is also
approximated as in Section 2. Let $N_{\e}$ be the total number of
discontinuities in the function $\overline{U}^{\e}$ and the
tangential angle functions of the wedge boundary. Choose a parameter
$\delta_{\e}>0$ controlling the maximum strength of rarefaction
fronts,  and $\hat\lambda$ (strictly larger than all the
characteristic speeds of (\ref{6.1.1})) that is the speed of
non-physical waves generated whenever the simplified method is used.
The strength of the non-physical waves is the error due to the
simplified Riemann solver.

{\bf A. Accurate Riemann solver}: The accurate Riemann solver is
just the solution to the Riemann problem (as in Section 2), except
every rarefaction wave $(w, R_i(w)(\alpha))$ is approximated by a
piecewise constant rarefaction fan.

{\bf B. Simplified Riemann solver}: For the weak waves, it is
exactly the same as in \cite{BJ}. When a weak wave interacts with
the large shock, the simplified Riemann solver is that we ignore the
strength of the weak wave, keep the strength of the strong shock,
and put the error in the non-physical
wave as follows:

Case 1 (A weak wave $(U_-, U_1)$ hits the large shock $(U_1, U_+)$
from below): We solve the Riemann problem $(U_-, U_+)$ in the
following way:
$$\left\{\begin{array}{ccc}
U_- \qquad &\mbox{for}&y/x<\Lambda(U_1, U_+),\\
U_1 \qquad &\mbox{for}&\Lambda(U_1, U_+)<y/x<\hat\lambda,\\
U_+ \qquad &\mbox{for}&y/x>\hat\lambda,
\end{array}\right.$$
where $\Lambda(U_1, U_+)$ is the speed of the strong shock;

Case 2 (A weak wave $(U_2, U_+)$ hits the large shock $(U_-, U_2)$
from above): We solve the Riemann problem $(U_-, U_+)$ in the
following way:
$$
\left\{\begin{array}{ccc}
U_- \qquad &\mbox{for}&y/x<\Lambda(U_-, U_2),\\
U_2 \qquad &\mbox{for}&\Lambda(U_-, U_2)<y/x<\hat\lambda,\\
U_+ \qquad &\mbox{for}&y/x>\hat\lambda,
\end{array}\right.$$
where $\Lambda(U_-, U_2)$ is the speed of the strong shock.

\subsection{The algorithm to construct the approximate solutions}
\setcounter{equation}{0} Given $\e$, we construct the approximate
solution $U^{\e}(x,y)$ as follows. When $x=0$, all the Riemann
problems in $\overline{U}^{\e}$ are solved by the accurate Riemann
solver. By slightly perturbing the speed of a wave, we can assume
that, at any time, we have at most one collision involving only two
incoming fronts. Let $\mu_\e$ be a fixed small parameter with
$\mu_\e\to 0$, as $\e\to 0$, which will be specified later. For
simplicity of notation, we will drop the index $i$ in $\alpha_i$ and
do not distinguish between $\alpha_i$ and $\alpha$ when there is no
ambiguity from now on; the same for $\beta$; also we use the same
notation $\alpha$ as a wave and its strength as before.

\noindent Case 1 (There is a collision between two weak waves with
strengths $\af$ and $\bt$ at some $x>0$, respectively): The Riemann
problem generated by this interaction is solved as follows:
\begin{itemize}
\item If $|\af\bt|>\mu_\e$ and the two waves are physical,
      then we use the accurate solver;
\item If $|\af\bt|<\mu_\e$ and the two waves are physical,
      or one wave is non-physical, then we use the simplified Riemann solver.
  \end{itemize}

\noindent Case 2 (There is a collision between the large shock and
one weak wave $\af$ at some $x>0$): The Reimann problem generated by
this interaction is solved as follows:
\begin{itemize}
\item If $|\af|>\mu_\e$ and the weak wave is physical, then we use the accurate solver;
\item If $|\af|<\mu_\e$ and the  weak wave is physical, or this wave is non-physical,
then we use the simplified Riemann solver.
\end{itemize}

\noindent Case 3 (The wave hits the boundary or the boundary
perturbs the flow): We use the accurate Riemann solver to solve the
lateral Riemann problem.

\subsection{Glimm's functional and interaction potential}
\setcounter{equation}{0} We now develop the Glimm-type functional
and interaction potential for the initial-boundary value problem by
carefully incorporating additional nonlinear waves generated from
the wedge boundary vortices.

\begin{definition}[Approaching waves]
{\rm (i)} We say that two weak fronts $\af$ and $\bt$, located at
points $x_{\af}<x_{\bt}$ and belonging to the characteristic
families $i_\af, i_\bt\in\{1, \dots, 4\}$ respectively, approach
each other if the following two conditions hold:
\begin{itemize}
\item $x_{\af}$ and $x_{\bt}$ both lay in one of the two intervals into
which {\bf R} is partitioned by the location of the large $1$-shock,
i.e. the waves both belong to $\Omega_-$ or $\Omega_+$;
\item Either $i_\af=i_\bt$ and one of them is a shock, or $i_\af>i_\bt$.
\end{itemize}
In this case we write $(\af, \bt)\in \mathscr{A}$.

{\rm (ii)} We say that a weak wave $\af$ of the characteristic
family $i_\af$ is approaching the large 1-shock if either $\af\in
\Omega_-$ and $i_\af\in\{1,2,3,4\}$, or  $\af\in \Omega_+$ and
$i_\af=1$.  We then write $\af\in \mathscr{A}_1$.

{\rm (iii)} We say that a weak wave $\af$ of the characteristic
family $i_\af$ is approaching the boundary if   $\af\in \Omega_+$
and $i_\af=4$.  We then write $\af\in \mathscr{A}_b$.
\end{definition}

For a weak wave $\af$ of $i$-family, we define its weighted strength
as
$$
b_\af=\left\{\begin{array}{ll}
\af \qquad \mbox{if}\;\af\in\Omega_+,\\
k_-\af \qquad \mbox{if}\;\af\in\Omega_-,
\end{array}\right.
$$
where $k_-=2\max\limits_{1\le i\le 4, 2\le j\le 4}\{K_{ij}\}$ for
coefficients $K_{ij}$ in Lemma \ref{prop6.3.4}.

\begin{definition}\label{Q}
The wave interaction potential $Q(x)$ is
\begin{eqnarray}
Q(x)&=&C^*\sum_{(\af, \bt)\in \mathscr{A}}|b_\af b_\bt| +
K^*\sum_{\af\in \mathscr{A}_1}|b_\af|+\sum_{\bt\in
\mathscr{A}_b}|b_\bt|
+\widetilde{K_{b0}}\sum_{a_k>x}|\o_k|\\
&=&Q_{\mathscr{A}}+Q_1+Q_b+Q_w,
\end{eqnarray}
where  $K^*\in (K_{s4},1)$ and $\widetilde{K_{b0}}>K_{b0}$.
\end{definition}

\begin{definition}
The total (weighted) strength of weak waves in $U^\e(x, \cdot)$ is
defined by
$$
V(x)=\sum_{\af}|b_\af|.
$$
The Glimm-type functional is defined by
\begin{equation}
\mathcal{F}(x)=V(x)+ \kappa Q(x)+|U^*(x)-U_0^+|+|U_*(x)-U_0^-|,
\end{equation}
where $\kappa>0$ is a constant to be specified later, the vectors
$U^*(x)$ and $U_*(x)$ are the above and below states of the large
shock respectively at ``time" $x$, and $U_0^+$ and $U_0^-$ are the
right and left states of the large shock at $x=0$, respectively.
\end{definition}

Note that $V$, $Q$, and $\mathcal{F}$ are constant between any pair
of subsequent interaction times. On the other hand, we can show
that, across an interaction ``time'' $x$, both $Q$ and $\mathcal{F}$
decrease.

\begin{proposition}\label{prop3.1}
If $TV (\bar{U}(\cdot))+TV(g'(\cdot))$ is sufficiently small, then,
for any $x>0$, $V(x)$ is sufficiently small and $TV (U^\e(x,
\cdot))$ is uniformly bounded in $\e>0$.
\end{proposition}
\begin{proof} Define
$$
\Delta\mathcal{F}(x)=\mathcal{F}(x^+)-\mathcal{F}(x^-),
$$
where $x^+$ and $x^-$ are the ``times'' right after and right before
the interaction ``time'', respectively.

Case 1 (Weak waves $\af$ and $\bt$ interact): Then $U^*(x)$ and
$U_*(x)$ do not change across this interaction time. Thus,
\begin{eqnarray*}
\mathcal{F}(x^+)-\mathcal{F}(x^-)&=&V(x^+)-V(x^-)+ \kappa
(Q(x^+)-Q(x^-))\\&\le &M_1|b_\af b_\bt|+ \kappa (-C^*|b_\af
b_\bt|+C^*|b_\af b_\bt| V(x^-)+ M_0|b_\af b_\bt|),
\end{eqnarray*}
for some constants $M_0$ and $M_1$ independent of $\e$.

Case 2 (Weak wave $\af$ of 4-family  interacts with the boundary):
\begin{eqnarray*}
\Delta\mathcal{F}(x)=K_{b4}\af -\af
+\kappa\left(C^*K_{b4}\,V(x^-)\,\af+K^*K_{b4}\af-\af\right).
\end{eqnarray*}

Case 3 (New 1-wave $\af$ produced by the boundary):
\begin{eqnarray*}
\Delta\mathcal{F}(x)=K_{b0}\o_k + \kappa\left(C^*K_{b0}\o_k\,
V(x^-)+K^*K_{b0}\o_k-\widetilde{K_{b0}}\o_k\right).
\end{eqnarray*}

The next two cases when $U^*(x)$ and $U_*(x)$ change across this
interaction ``time''. Then

Case 4 (Weak wave $\af$ of $i$-family interacts with the strong
shock from below):
\begin{eqnarray*}
\Delta\mathcal{F}(x)&=&V(x^+)-V(x^-)+|U^*(x^+)-U^*(x^-)|\\
&&
+|U_*(x^+)-U_*(x^+)|+K(Q(x^+)-Q(x^-))\\
&=&\sum_{j=1,2,3,4}K_{ji}\ep_\af
-b_\af+K\big(C^*\sum_{j=2,3,4}K_{ji}\,
V(x^-)\,\af-b_{\af}+K_{4i}\af\big);
\end{eqnarray*}

Case 5 (Weak wave $\af$ of 1-family  interacts with the strong shock
from above):
\begin{eqnarray*}
\Delta\mathcal{F}(x)&=&\sum_{j=1,2,3,4}K_{si}\af -\af
+K\big(C^*\sum_{j=2,3,4}K_{sj}\,V(x^-)\,\af
-K^*{\af}+K_{s4}\af\big).
\end{eqnarray*}
In these cases, $K_{s4}<K^*<1$, $b_\af \ge 2\max \{K_{ji}\}|\af|$
due to the choice of the weight $k_-$, and $C^*>M_0>0$ is a constant
that is not small.

We now prove
$$
V(x)\ll 1 \qquad\text{for all } \, x>0.
$$

Case 1 ($x>0$ is the first interaction ``time'' $x_1$): Since
$V(x^-_1)=V(0)\le TV(\bar{U}(\cdot))\ll 1$ and $\sum_{x=0}^\infty
\o_k\le TV(g'(\cdot))\ll 1$ for all Cases 1--5, we find that, when
$\kappa$ is larger enough and $\mu_\e$ is sufficiently small,
$$
\Delta\mathcal{F}(x_1)\le 0,\qquad i.e.
\qquad\mathcal{F}(x_1^+)\le\mathcal{F}(x_1^-)=\mathcal{F}(0).
$$
Therefore,
\begin{eqnarray*}
V(x^+_1)&\le& \mathcal{F}(x_1^+)\le \mathcal{F}(0)\le V(0)+ \kappa Q(0)\\
&=&V(0)+\kappa \big(C^*V^2(0)+V(0)+\widetilde{K_{b0}}\sum_{x=0}^\infty \o_k\big)\\
&\le&C\big(V(0)+\sum_{x=0}^\infty \o_k\big)\ll 1.
\end{eqnarray*}

Case 2 ($V(x^-_{k})\ll 1$ and
$\mathcal{F}(x_k^+)\le\mathcal{F}(x_k^-)$ for any $k<n$): Then, for
the next interaction ``time'' $x_n$, similarly to Case 1, we also
have
$$
\Delta\mathcal{F}(x_n)\le 0,\qquad i.e.
\qquad\mathcal{F}(x_n^+)\le\mathcal{F}(x_n^-)=\mathcal{F}(x_{n-1}^+).
$$
Thus, we have
\begin{eqnarray*}
&&V(x^+_n)+|U^*(x^+_n)-U_0^+|+|U_*(x^+_n)-U_0^-|\\
&&\le \mathcal{F}(x_n^+)\le\mathcal{F}(x_n^-)=\mathcal{F}(x_{n-1}^+)
\le\ldots\le\mathcal{F}(0)
=V(0)+ \kappa Q(0)\\
&&=V(0)+\kappa \big(C^* V^2(0)+V(0)+\widetilde{K_{b0}}\sum_{x=0}^\infty \o_k\big)\\
&&\le C\big(V(0)+\sum_{x=0}^\infty \o_k\big)\ll 1.
\end{eqnarray*}
Therefore, $V(x)\ll 1$ is proved for all $x$, since $C$ is
independent of $x$. Then
\begin{equation}
TV \{U(x,\cdot)\}\approx
V(x)+|U^*(x)-U_0^+|+|U_*(x)-U_0^-|+|\S_0|=O(1).
\end{equation}
\end{proof}

\begin{lemma}\label{lemma3.2} For any sufficiently small fixed $\e>0$,
the number of wave fronts in $U^\e(x,y)$ is finite so that the
approximate solutions $U^\e(x,y)$ are defined for all $x$.
\end{lemma}

\begin{proof} Recall that the total interaction potential $Q(x)$ is constant
except decreasing when it crosses an interaction time. From Cases
1--5 in Proposition \ref{prop3.1}, we have known that $V(t)\ll 1$.
Therefore, we can find some $c\in (0,1)$ so that
\begin{equation}\label{Delta Q(t)}
\begin{array}{ll}
&\Delta Q(x)=Q(x^+)-Q(x^-)\\
&\qquad\quad\le \left\{\begin{array}{ll}
-c|b_\af b_\bt|\qquad &\mbox{if both waves $\af$ and $\bt$ are weak,}\\
-c|b_\af|\qquad &\mbox{if weak wave $\af$ hits the strong shock or the boundary, }\\
-c|\o_k|\qquad &\mbox{if the angle of the boundary changes.}
\end{array}\right.
\end{array}
\end{equation}
The following argument is similar to that in \cite{BJ}: $Q$
decreases for each case and $Q(0)$ is bounded; When the interaction
potential between the incoming waves is greater than $\mu_\e$, $Q$
decreases by at least $c\mu_\e$ in these interactions, by the bound
in (\ref{Delta Q(t)}); Following the wave front tracking, new
physical waves can be only generated by this kind of interactions,
which implies that the number of the waves is finite; Since
non-physical waves are produced only when physical waves interact,
the number of non-physical waves is also finite; and, since two
waves can interact only once, the number of interactions is also
finite. Therefore, the approximate solutions are defined for all
$x>0$.
\end{proof}

Similar to \cite{BJ}, we have following lemma.
\begin{lemma}
The total strength of all non-physical waves at any $x$ is of the
order $O(1)(\delta_\e+\mu_\e)$.
\end{lemma}

Following the framework of the wave front tracking in \cite{Br,BJ}
and Lemmas \ref{prop3.1}--\ref{lemma3.2}, we obtain the existence of
global entropy solutions to (\ref{eq1.1}) and
(\ref{IC})--(\ref{BC}).

\begin{theorem}
If $TV (\bar{U}(\cdot))+TV(g'(\cdot))$ is sufficiently small, then
there exists a global entropy solution in BV of problem {\rm
(\ref{eq1.1})} and {\rm (\ref{IC})}--{\rm (\ref{BC})} of
initial-boundary value type in the sense of Definition {\rm 2.1}.
\end{theorem}

\section{The Lyapunov functional}
\setcounter{equation}{0} We now follow the approach of
\cite{BLY,LT,LY} to construct the Lyapunov functional $\Phi(U, V)$
by incorporating additional new waves generated from the wedge
boundary vortices, which is equivalent to the $L^1$-distance:
$$
C_1^{-1}\|U(x,\cdot)-V(x,\cdot)\|_{L^1}\le\Phi(U, V)\le
C_2\|U(x,\cdot)-V(x,\cdot)\|_{L^1},
$$
$$
\Phi(U(x_2,\cdot), V(x_2,\cdot))-\Phi(U(x_1,\cdot), V(x_1,\cdot))
\le C_3\ep (x_2-x_1)\qquad\mbox{for any}\,\, x_2>x_1>0,
$$
for some constants $C_i, i=1,2,3$, where $U$ and $V$ are two
approximate solutions obtained by the wave front tracking, and the
small parameter $\ep$ controls the following three types of errors:
\begin{itemize}
\item Errors in the approximation of initial data and boundary;
\item Errors in the speeds of shock and rarefaction fronts;
\item The maximum strength of rarefaction fronts;
\item The total strength of all non-physical waves.
\end{itemize}
When $x$ is fixed, for each $y$, the connection $U(y)$ with $V(y)$
always moves along Hugoniot curves $S_1,C_2, C_3,$ and $S_4$ in the
phase space. We call $p_i(y)$ the strength of the $i$-th
discontinuity wave, which is determined by $U(y)$ and $V(y)$ as
follows:
\begin{itemize}
\item If $U(y)$ and $V(y)$ are both in $\Omega_-$ or in $\Omega_+$,
     then start from $U(y)$ moving along Hugoniot curves and end at $V(y)$;
\item If $U(y)$ is in $\Omega_-$, $V(y)$ is in $\Omega_+$,
      then also start from $U(y)$ moving along Hugoniot curves and end at $V(y)$.
\item If $V(y)$ is in $\Omega_-$ and $U(y)$ is in $\Omega_+$, then start from $V(y)$
      moving along Hugoniot curves and end at $U(y)$.
\end{itemize}
Now we define the weighted $L^1$ strength:
\begin{equation}\label{q(x)}
q_i(y)=\left\{\begin{array}{ll} c_i^bp_i(y)\qquad&\mbox{if  $U(y)$
and $V(y)$ are both in $\Omega_-$},\\
c_i^mp_i(y)\qquad&\mbox{if
$U(y)$ and $V(y)$ are in different domains},\\
c_i^ap_i(y)\qquad&\mbox{if  $U(y)$ and $V(y)$ are both in
$\Omega_+$},
\end{array}\right.
\end{equation}
where the constants $c_i^b$, $c_i^m$, and $c_i^a$ are to be
determined later. Then we define the Lyapunov functional:
\begin{equation}
\Phi(U, V)=\sum_{i=1}^4\int_{-\infty}^{g(x)} |q_i(y)|W_i(y)dy,
\end{equation}
with
\begin{equation}
W_i(y)=1+\kappa_1A_i(y)+\kappa_2(Q(U)+Q(V)).
\end{equation}
Here $\kappa_1$ and $\kappa_2$ are two constants to be defined
later, $Q$ is the total wave interaction potential defined in
Definition \ref{Q}, $A_i(y)$ is the total strength of waves in $U$
and $V$ which approach the $i$-wave $q_i(y)$ defined by
\begin{equation}
A_i(y)=B_i(y)+D_i(y)
+\left\{\begin{array}{ll}
C_i(y)\quad\mbox{if} \;q_i(y)\; \mbox{is small}, \\
F_i(y)\quad\mbox{if}\; i=1 \;\mbox{and}\;  q_1(y)=B\; \mbox{is
large},
\end{array}\right.
\end{equation}
where the ``small" or ``large" describes the waves that connect the
states in the same or in the distinct domains $\Omega^-$ and
$\Omega^+$, respectively, and
\begin{eqnarray*}
B_i(y)&=&\bigg(\sum_{{\af\in\mathcal{J}(U)\cup\mathcal{J}(V)}\atop{y_\af<y,
i<k_\af\le 4}}
+\sum_{{\af\in\mathcal{J}(U)\cup\mathcal{J}(V)}\atop {y_\af>y, 1\le k_\af< i}}\bigg)
|\af|,\\
C_i(y)&=&\left\{\begin{array}{cc}
 \left(\sum_{{\af\in\mathcal{J}(U)}\atop {y_\af<y, k_\af=i}}
    +\sum_{{\af\in\mathcal{J}(V)}\atop{y_\af>y, k_\af=i}}\right)|\af|
     \qquad\mbox{if}\;q_i(y)<0,\\
 \left(\sum_{{\af\in\mathcal{J}(V)}\atop{y_\af<y, k_\af=i}}
     +\sum_{{\af\in\mathcal{J}(U)}\atop {y_\af>y, k_\af=i}}\right)
     |\af|\qquad\mbox{if}\;q_i(y)>0,\\
\end{array}\right. \\
F_i(y)&=&\bigg(\sum_{{{\af\in\mathcal{J}(U)\cup\mathcal{J}(V)}\atop{
y_\af<y, k_\af=1}}\atop{ {\mbox{\tiny both states joined by $\af$}
}\atop {\mbox{\tiny are located in
$\Omega_-$}}}}+\sum_{{{\af\in\mathcal{J}(U)\cup\mathcal{J}(V)}\atop{
y_\af>y, k_\af=1}}\atop{ {\mbox{\tiny both states joined by $\af$}
}\atop {\mbox{\tiny are located in $\Omega_+$}}}}\bigg)|\af|.
\end{eqnarray*}\\
For fixed $x$, $\mathcal{J}=\mathcal{J}(U)\cup\mathcal{J}(V)$ is the
set of all weak waves in $U$ and $V$,  $\af$ is the strength of wave
$\af\in \mathcal{J}$, located at point $y_\af$ and belonging to the
characteristic family $k_\af$.
\\
\\
\begin{tabular}{l|l|l|l}
$D_i(y)=$& $U$, $V$ are both in $\Omega_-$& $U$, $V$ are  in different domains&$U$,
$V$ are both in $\Omega_+$\\
\hline
$D_1(y)$&$B$&0&$B$\\
$D_{2,3}(y)$&$B$&$B$&0\\
$D_4(y)$&$B$&$B$&$B$
\end{tabular}
\\
\\
Since, for any $U(x,\cdot)$, $V(x, \cdot)\in BV\cap L^1$ and
$TV(\bar{U}(\cdot))+TV(\bar{V}(\cdot))+ TV(g'(\cdot))$ is
sufficiently small, we have
\begin{eqnarray*}
&&C^{-1}_0\|U(x,\cdot)-V(x,\cdot)\|_{L^1}
\le\sum_{i=1}^4\int_{-\infty}^{g(x)}
|q_i(y)|dy\le C_0\|U(x,\cdot)-V(x, \cdot)\|_{L^1},\\
&&1\le W_i(y)\le C_0,\qquad i=1,2,3,4,
\end{eqnarray*}
for some constant $C_0$ independent of $x$ and $\e$. Therefore, for
any $x\ge 0$,
\begin{equation}\label{Phieq}
C_1^{-1}\|U(x,\cdot)-V(x, \cdot)\|_{L^1}\le \Phi(U, V)\le
C_2\|U(x,\cdot)-V(x,\cdot)\|_{L^1},
\end{equation}
where $C_1$ and $C_2$ depend only on
$TV(\bar{U}(\cdot)+TV(\bar{V}(\cdot)+TV (g'(\cdot))$ and the
strength of the strong shock, which are independent of $x$.

Now we examine how the Lyapunov functional $\Phi$ evolves in the
flow direction $x>0$. Denote $\lambda_i$ the speed of the $i-$wave
$q_i(x)$ (along the Hugoniot curve in the phase space). At a time
$x$ which is not the interaction time of the waves either in $U$ or
$V$,
\begin{eqnarray*}
&&\frac{d}{dx}\Phi(U(x), V(x))\\
&&=\sum_{\af\in\mathcal{J}}\sum_{i=1}^4
\big(|q_i(y_\af-)|W_i(y_\af-)-|q_i(y_\af+)|W_i(y_\af+)\big)\dot
y_\af
+\sum_{i=1}^4|q_i(b)|W_i( b)\dot y_b\\
&&=\sum_{\af\in\mathcal{J}}\sum_{i=1}^4
\big(|q_i(y_\af-)|W_i(y_\af-)(\dot y_\af-\lambda_i(y_\af-))
-|q_i(y_\af+)|W_i(y_\af+)(\dot y_\af-\lambda_i(y_\af+))\big)\\
&&\quad +\sum_{i=1}^4|q_i(b)|W_i( b)(\dot y_b-\lambda_i(b)),
\end{eqnarray*}
where $\dot y_\af$ is the speed of the discontinuity at wave
$\af\in\mathcal{J}$, $b=g(x)-$ stands for the points on the
boundary, and $\dot y_b$ is the slope of the boundary. Define
\begin{eqnarray}
E_{\af, i}&=&|q_i^+|W_i^+(\lambda_i^+-\dot y_\af)-|q_i^-|W_i^-(\lambda_i^--\dot y_\af),\\
E_{b, i}&=&|q_i(b)|W_i( b)(\dot y_b-\lambda_i(b)),
\end{eqnarray}
where $q_i^\pm=q_i(y_\af\pm)$, $W_i^\pm=W_i(y_\af\pm)$, and
$\lambda_i^\pm=\lambda_i(y_\af\pm)$. Then
\begin{equation}
\frac{d}{dx}\Phi(U(x), V(x))=\sum_{\af\in\mathcal{J}}\sum_{i=1}^4
E_{\af, i}+\sum_{i=1}^4E_{b, i}.
 \end{equation}
Our main goal is to establish the following bounds:
\begin{eqnarray}
&&\sum_{i=1}^4E_{b, i}\le 0 \qquad\mbox{near the boundary,}\label{1}\\
&&\sum_{i=1}^4 E_{\af, i}\le 0\qquad
  \mbox{when $\af$ is a strong shock wave in $\mathcal{J}$},\label{2}\\
&&\sum_{i=1}^4 E_{\af, i}\le O(1)|\af|
  \qquad\mbox{when $\af$ is a non-physical wave
                    in $\mathcal{J}$},\label{3}\\
&&\sum_{i=1}^4 E_{\af, i} \le O(1)\varepsilon|\af|\qquad \mbox{when
$\af$ is a weak wave in $\mathcal{J}$}.\label{4}
\end{eqnarray}
{}From (\ref{1})--(\ref{4}), we have
 \begin{equation}
 \frac{d}{dx}\Phi(U(x), V(x))\le O(1)\ep.
 \label{Phi}
 \end{equation}
If the constant $\kappa_2$ in the Lyapunov functional is chosen
large enough, by the Glimm interaction estimates, all weight
functions $W_i(y)$ decrease at each time where two fronts of $U$ or
two fronts of $V$ interact.  By the self-similar property of the
Riemann solutions, $\Phi$ decreases at this time. Integrating
(\ref{Phi}) over interval $[0, x]$, we obtain
\begin{equation}\label{PhiUV}
\Phi(U(x), V(x))\le \Phi(U(0), V(0))+O(1)\ep \,x.
\end{equation}
In Section 5, we prove (\ref{1})--(\ref{4}).

\section {Estimates for the $L^1$ Stability}
\setcounter{equation}{0}

For the case that the weak wave $\af\in
\mathcal{J}:=\mathcal{J}(U)\cup\mathcal{J}(V)$ appears when $U$ and
$V$ both in $\Omega_-$ or $\Omega_+$  and for the case of   the
non-physical waves in  $\mathcal{J}$, estimates (\ref{3})--(\ref{4})
can be obtained if $|B/\S_0|$ is small enough and $\kappa_1$ is
large enough, by following \cite{BLY}.

Therefore, in this section, we focus the other cases. Cases 1--3
below are all related to the strong shock and depend on the wave
jump $\af$ in $U$ or $V$; and, by carefully adjusting the
coefficients $c_i$ and especially relying on estimate \eqref{key},
we can obtain desirable results, which is similar to the Cauchy
problem discussed in \cite{LT}. Case 4 is the case near the
boundary, which is different from those for the Cauchy problem.

Case 1 (Cross the first strong shock $\af$ in $U$ or $V$): For this
case,
\begin{eqnarray*}
E_1&=&B W_1^+(\lambda_1^+-\dot y_\af)-|q_1^-|W_1^-(\lambda_1^--\dot y_\af)\\
&\le&O(1) B\sum_{i=1}^4|q_i^-|-\kappa_1B|q_1^-|\,|\lambda_1^--\dot
y_\af|,
\end{eqnarray*}
and
\begin{eqnarray*}
\sum_{i=2}^4E_i&=&\sum_{i=2}^4\big(|q_i^-|(\lambda_i^--\dot
y_\af)(W_i^+- W_i^-)
 + W_i^+(|q_i^+|(\lambda_i^+-\dot y_\af)-|q_i^-|(\lambda_i^--\dot y_\af))\big)\\
&\le&\sum_{i=2}^4\kappa_1 B|q_i^+||\lambda_i^+-\dot y_\af|-
\frac{3}{4}\sum_{i=2}^4 \kappa_1 B|q_i^-||\lambda_i^--\dot y_\af|.
\end{eqnarray*}
Therefore, when $\kappa_1$ is large enough, we have
$$
\sum_{i=1}^4E_i\le\sum_{i=2}^4\kappa_1 B|q_i^+||\lambda_i^+-\dot
y_\af|-\sum_{i=1}^4\frac{1}{2}\kappa_1 B|q_i^-||\lambda_i^--\dot
y_\af|.
$$
In Lemma \ref{prop6.3.4}, let $\af_i=p_i^-$ and $\bt_i=0$, which
implies $\delta_i=p_i^+$. We know that
$$
p_k^+ \approx O(1)\sum_{i=1}^4|p_i^-|, \qquad k=2,3,4.
$$
For the weighted $L^1$ strength $q_i(y)$ in (\ref{q(x)}), when
$c_i^b, 1\le i\le 4,$ are larger enough relatively to $c_i^m,
i=2,3,4$, we can obtain (\ref{2}).

Case 2 (Cross the weak wave $\af$ in between the two strong shocks
in $U$ and $V$): For this case,
\begin{eqnarray*}
E_1&=&B\left((W_1^+-W_1^-)(\lambda_1^\pm-\dot y_\af)
  +W_1^\mp(\lambda_1^\pm-\lambda_1^\mp)\right)\\
&\le&B\left(-\kappa_1|\af||\lambda_1^+-\dot y_\af|+O(1)|\af|\right).
\end{eqnarray*}
For $i=2,3,4$,
\begin{eqnarray*}
E_i&=&|q_i^\pm|(W_i^+-W_i^-)(\lambda_i^\pm-\dot y_\af)
    +W_i^\mp\big(|q_i^+|(\lambda_i^+-\dot y_\af)-|q_i^-|(\lambda_i^--\dot y_\af)\big)\\
&\le&\kappa_1|q_i^\pm||\af||\lambda_i^\pm-\dot y_\af|
    + \kappa_1B\big((|q_i^+|-|q_i^-|)(\lambda_i^+-\dot y_\af)
      +|q_i^-|(\lambda_i^+-\lambda_i^-)\big)\\
&\le&\kappa_1|q_i^\pm||\ep_\af||\lambda_i^\pm-\dot y_\af|
   + \kappa_1B\big((|q_i^+|-|q_i^-|)(\lambda_i^+-\dot
   y_\af)+O(1)|q_i^-||\af|\big).
\end{eqnarray*}
Then we have
$$
\sum_{i=1}^4E_i\le \kappa_1
O(1)\big(-|\af|+|\af|\sum_{k\neq1}(|q_k^+|+|q_k^-|)
 +\sum_{k\neq1}(|q_k^+|-|q_k^-|)\big)+O(1)|\af|.
$$
Since $||q_k^+|-|q_k^-||\le|q_k^+-q_k^-|\le O(1)|\af|
\;\mbox{when}\;k\neq 1$, we can obtain $\sum_{i=1}^4E_i\le 0$ if all
the weights $c_i^m$ are sufficiently small and $\kappa_1$ is large
enough.

Case 3 (Cross the second strong shock $\af$ in $U$ or $V$): For this
case, by Lemma \ref{prop6.3.3}, we have
\begin{equation}\label{p4}
p_4^-=p_4^++\widetilde {K_{s4}}p_1^+.
\end{equation}
Since we have (\ref{key}),  the following lemma can be easily
obtained.
\begin{lemma}\label{5.1}
There exist $c_1^a$, $c_4^a$, and  $\gamma$ such that
\begin{eqnarray}
&&\frac{c_1^a}{c_4^a}<1,\qquad\qquad\label{con1}\\
&&\frac{c_4^a}{c_1^a}\widetilde{K_{s4}}\frac{|\lambda_4^+-\S|}{|\lambda_1^+-\S|}<\gamma<1.
\label{con2}
\end{eqnarray}
\end{lemma}

With Lemma \ref{p4}, then we estimate $E_i$:
\begin{eqnarray*}
E_1&=&-B W_1^-(\lambda_1^--\dot y_\af)+|q_1^+|W_1^+(\lambda_1^+-\dot y_\af)\\
&\le&O(1) B|q_1^+|-\kappa_1B|q_1^+|\, |\lambda_1^+-\dot y_\af|\\
&=&O(1) B|q_1^+|-\kappa_1Bc_1^a|p_1^+|\, |\lambda_1^+-\dot y_\af|,
\end{eqnarray*}
and, for $i=2,3,$
\begin{eqnarray*}
E_{i} &=&|q_i^-|(\lambda_i^--\dot y_\af)(W_i^+-W_i^-)
 +W_i^+\big(|q_i^+|(\lambda_i^+-\dot y_\af)-|q_i^-|(\lambda_i^--\dot y_\af)\big)\\
&\le&-\kappa_1B|q_i^-|(\lambda_i^--\dot y_\af)+O(1)|q_i^+|\\
&\le&-\kappa_1B|q_i^-|(\lambda_i^--\dot
y_\af)+O(1)(|q_i^-|+|q_1^+|).
\end{eqnarray*}
By (\ref{p4}) and (\ref{con2}),
 \begin{eqnarray*}
E_4 &=&|q_4^-|(\lambda_4^--\dot y_\af)(W_4^+-W_4^-)
+W_4^+\big(|q_4^+|(\lambda_4^+-\dot y_\af)-|q_4^-|(\lambda_4^--\dot y_\af)\big)\\
&\le&\kappa_1B\left(c_4^a|p_4^-|(\lambda_4^+-\dot y_\af)
 +c_4^a\widetilde{K_{s4}}|p_1^+|(\lambda_4^+-\dot y_\af)-c_4^m|p_4^-|(\lambda_4^--\dot y_\af)\right)\\
&\le&\kappa_1B\left(c_4^a|p_4^-|(\lambda_4^+-\dot y_\af)+\gamma
c_1^a|p_1^+|\,|\lambda_1^+-\dot y_\af|-c_4^m|p_4^-|(\lambda_4^--\dot
y_\af)\right).
\end{eqnarray*}
From above, if we choose $c_4^a$ is small enough relatively to
$c_4^m$ and choose $k_1$ is large enough, then we obtain
\begin{eqnarray*}
\sum_{i=1}^4 E_i&=&-(1-\gamma)\kappa_1B|q_1^+|\,|\lambda_1^+-\dot y_\af|+O(1)\, |q_1^+|\\
&&+ \kappa_1B(c_4^a|p_4^-|(\lambda_4^+-\dot
y_\af)-c_4^m|p_4^-|(\lambda_4^--\dot y_\af))\\
&&+\sum_{i=2}^3(-\kappa_1B|q_i^-|(\lambda_i^--\dot
y_\af)+O(1)\cdot|q_i^-|)\le 0.
\end{eqnarray*}

Case 4 (near the boundary): For the previous cases, all the desire
results depend on the cancellation between the two sides of a wave
in $\mathcal{J}$.  However, it is not the case near the boundary
since there is only one side near the boundary. Then we exploit the
exclusive property of the boundary condition (\ref{BC}): the flows
of $U$ and $V$ are tangent to the boundary, which implies that they
must be parallel to each other. Then we solve the Riemann problem
determined by $U(b)$ and $V(b)$.

\begin{proposition}\label{Eb}
Let $U(b)=(\bar{u}, \bar{v}, \bar{p},\bar{\rho})$ and
$V(b)=(\hat{u}, \hat{v}, \hat{p}, \hat{\rho})$ be both in $O(U_+)$,
$v_1/u_1= v_2/u_2=\dot y_b$, and $\bar{v}, \hat{v}\approx 0$. Let
$p_i(b)$ and $\lambda_i$ be the strength and speed of the $ith$
shock in the Riemann problem determined by $U(b)$ and $V(b)$. Then
\begin{eqnarray}
&&|\dot y_b-\lambda_i| \sim |p_1(b)|,\qquad i=2, 3,\\
&&|p_1(b)|\le |p_4(b)|+ O(1)|p_2(b)||\lambda_2-\dot y_b|+|p_1(b)| O(1)|\dot y_b|,\\
&&|p_4(b)|= O(1)|p_1(b)|.
\end{eqnarray}
\end{proposition}

\begin{proof} We divide the proof into two cases.

Case 1. $p_1(b)=0=p_4(b)$ which corresponds the case
$\bar{p}=\hat{p}$: Starting from $U_b$, go along the curves of the
second and third families to reach $V_b$. These two families are the
contact Hugoniot curves,
and $\lambda_{2}$ and $\lambda_3$ are constant along the curves.
Since $\lambda_{2,3}=v/u$, ${\bf r}_2=(1, v/u, 0, 0)^\top$, and
${\bf r}_3=(0, 0, 0, \rho)^\top$, $v/u$ keeps unchanged as the
initial value $v(U_b)/u(U_b)$, i.e. $\dot y_b$ in this process.
Therefore, $\lambda_{2,3}=\dot y_b$, i.e.,
$$
\dot y_b-\lambda_{2,3}=0.
$$

Case 2. $p_1(b)\neq 0$ which corresponds to $\bar{p}_1\ne
\hat{p}_4$. Starting from $U(b)$, go along the first Hugoniot curve
to reach $U_1$, then possibly along the second curve to reach $U_2$,
the third curve to reach $U_3$, and the fourth Hugoniot curve to
reach $V(b)$.

We project $(u,v, p, \rho)$ onto the $u-v$ plane to see the relation
among $p_1(b)$, $p_2(b)$, $p_3(b)$, and $p_4(b)$ more clearly.
Denote ${\bf r}_1|_u$ the projection of ${\bf r}_1$ onto the $u$
axis, ${\bf r}_2|_{(u,v)}$ the projection of ${\bf r}_2$ onto the
$u-v$ plane; and the others are defined similarly. At $U_+$, we have
$$
{\bf r}_1|_u=-{\bf r}_4|_u, \,\, {\bf r}_1|_v={\bf r}_4|_v,\,\, {\bf
r}_1|_{(p,\rho)}=-{\bf r}_4|_{(p,\rho)}, \,\,{\bf r}_2={\bf
r}_2|_{(u,v)},\,\, {\bf r}_3|_{(u,v)}=0.
$$
The first observation is $p_4(b)\neq 0$. Since ${\bf
r}_1|_{(u,v)}=k_1(-\lambda_1, 1)^\top$, the characteristic speed is
finite and $\dot y_b\approx 0$, so we always have
$\frac{-1}{\lambda_1}>\dot y_b$ near $U_+$, i.e. the derivative
$dv/du$ along the 1st curve is always larger than $\dot y_b$ in the
$u-v$ plane.  Therefore, $v(U_1)/u(U_1)\neq v(U_b)/u(U_b)$. On the
other hand, we have  $v(U_1)/u(U_1)=v(U_2)/u(U_2)=v(U_3)/u(U_3)$ and
$v(V_b)/u(V_b)=v(U_b)/u(U_b)$. Therefore,
$v(U_1)/u(U_1)=v(U_2)/u(U_2)=v(U_3)/u(U_3)\neq v(V_b)/u(V_b)$. To
reach $V_b$, there must be some distance along the 4th Hugoniot
curve. Thus, $p_4\neq 0$.

On the $u-v$ plane, we define the signed length  of
$(U_1-U_b)|_{(u,v)}$ and $(V_b-U_3)|_{(u,v)}$  by $l_1$ and $l_4$ as
follows:
$$
l_1=\left\{\begin{array}{ll} \|(U_1-U_b)|_{(u,v)}\| \qquad &\mbox{if
} p_1>0,\\
-\|(U_1-U_b)|_{(u,v)}\|\qquad &\mbox{if } p_1<0;
\end{array}\right.
$$
and
$$
l_4=\left\{\begin{array}{ll} \|(V_b-U_3)|_{(u,v)}\|\qquad&\mbox{if }
p_4>0, \\
-\|(V_b-U_3)|_{(u,v)}\|\qquad&\mbox{if } p_4<0.
\end{array}\right.
$$

The second observation is
$$
|\lambda_2-\dot y_b|=O(1)|l_1|=O(1) |p_1(b)|,
$$
and, since $\lambda_2=v(U_1)/u(U_1)=v(U_2)/u(U_2)=\lambda_3$, we
also have
$$
|\lambda_3-\dot y_b|=O(1) |p_1(b)|.
$$

The third observation is
$$
-l_4=p_2(b)\cdot O(1)(\lambda_2-\dot y_b)+\widetilde l,
$$
where $\widetilde l\, cos \theta_1=l_1 cos \theta_2$, $\theta_1$ is
the angle between $(1, \dot y_b)$  and ${\bf r}_4|_{(u,v)}$,
$\theta_2$ is the angle between ${\bf r}_1|_{(u,v)}$ and $(1, \dot
y_b)$, $\theta_1=\theta_2+2\bt$ for $\bt=arctan\, \dot y_b$, and
\begin{eqnarray*}\widetilde l&=&l_1 \frac{cos \theta_2}{cos \theta_1}
=l_1 \frac{cos (\theta_1-2\bt)}{cos \theta_1}
=l_1 \frac{cos \theta_1 cos(2\bt)+sin \theta_1 sin(2\bt)}{cos \theta_1}\\
&=&l_1\big(cos(2\bt)+O(1)sin(2\bt)\big)=l_1\big(1+O(1)\bt\big)
=l_1\big(1+O(1)\dot y_b\big).
\end{eqnarray*}
Therefore, we have
$$
-l_4=O(1) p_2(b)(\lambda_2-\dot y_b)+l_1(1+O(1)\dot y_b).
$$
Since ${\bf r}_1|_{(u,p,\,\rho)}=-{\bf r}_4|_{(u,p,\,\rho)}$ and
${\bf r}_1|_v={\bf r}_4|_v$  at $U_+$, we have
$$
\frac{l_1}{p_1}=\frac{l_4}{p_4}.
$$
Then we obtain
\begin{equation}\label{p_4}
-p_4(b)=O(1) p_2(b)(\lambda_2-\dot y_b)+p_1(b)(1+O(1)\dot y_b).
\end{equation}
Therefore, from (\ref{p_4}), we obtain
$$
|p_1(b)|\le |p_4(b)|+ O(1)|p_2(b)||\lambda_2-\dot y_b|+|p_1(b)|
O(1)|\dot y_b|,
$$
and
$$
|p_4(b)|=O(1)|p_1(b)|.
$$
\end{proof}

By Proposition \ref{Eb},
\begin{eqnarray*}
E_{b, 1}&=&|q_1(b)|W_1( b)(\dot y_b-\lambda_1)
=c_1^a|p_1(b)|\kappa_1B(-\lambda_1)+O(1)|p_1(b)|\\&\le&c_1^a|p_4(b)|\kappa_1B(-\lambda_1)+O(1)|p_1(b)|,\\
E_{b, i}&=&|q_i(b)|W_i( b)(\dot
y_b-\lambda_i)=c_i^a|p_i(b)|O(1)(\dot y_b-\lambda_i)
=O(1)|p_1(b)|,\qquad i=2,3,\\
E_{b, 4}&=&|q_4(b)|W_4( b)(\dot y_b-\lambda_4)
=c_4^a|p_4(b)|\kappa_1B\,\lambda_1+O(1)|p_4(b)|\\
&=&c_4^a|p_4(b)|\kappa_1B\,\lambda_1+O(1)|p_1(b)|.
\end{eqnarray*}
From Lemma \ref {5.1}, we can find $c_1^a$ and $c_4^a$ such that
$c_1^a<c_4^a$. Then, when $\kappa_1$ is large enough, we conclude
\begin{eqnarray*}\sum_{i=1}^4E_{b,
i}&=&(c_1^a-c_4^a)\,|p_4(b)|\,
\kappa_1B(-\lambda_1)+O(1)|p_1(b)|\\&\le& (c_1^a-c_4^a)\,
O(1)|p_1(b)|\, \kappa_1B(-\lambda_1)+O(1)|p_1(b)|\le
0.\end{eqnarray*}

\section{Semigroup}
\setcounter{equation}{0} We now prove the existence of the semigroup
generated by the wave front tracking method.

\begin{proposition} If $TV(\bar{U}(\cdot))+ TV(g'(\cdot))$ is sufficiently small,
then, the map $(\bar{U}(\cdot), x)\mapsto U^\ep(x,\cdot):=S^\ep
_x(\bar{U}(\cdot))$ produced by the wave front tracking method is a
uniformly Lipschitz continuous semigroup with following properties:
\begin{enumerate}
\item[(i)]  $S^\ep_0 \bar{U}=\bar{U} $,  $S^\ep_{x_1}S^\ep_{x_2}\bar{U}=S^\ep_{x_1+x_2}\bar{U}$;
\item[(ii)] $||S^\ep_x \bar{U}-S^\ep_x \bar{V}||_{L^1}\leq C||\bar{U}-\bar{V}||_{L^1}+C\ep\, x$.
\end{enumerate}
\end{proposition}

\begin{proof}
Property (i) is obvious since $S^\ep$ is produced by the wave front
tracking method. Then we see property (ii).

Let $\{U^\ep\}$ and $\{V^\ep\}$ be the front tracking
$\ep$-approximate solutions of (\ref{eq1.1}) and
(\ref{IC})--(\ref{BC}) with initial data functions $\bar{U}(\cdot)$
and $\bar{V}(\cdot)$, respectively. By (\ref{PhiUV}) and
(\ref{Phieq}), we obtain that, for any $x\ge 0$,
\begin{eqnarray*}
\|U^\ep(x)-V^\ep(x)\|_{L^1}
&\le& C_1 \Phi(U^\ep(x), V^\ep(x))\\
 &\le &C_1 \Phi(U^\ep(0), V^\ep(0))+C_1\, O(1)\ep x\\
 &\le & C_1C_2\|\bar{U}-\bar{V}\|_{L^1}+ C_1\, O(1)\ep x.
 \end{eqnarray*}
This establishes the Lipschitz continuity of the $\ep$-semigroup.
\end{proof}

\begin{definition} Given $\delta_0>0$, define the domain
${\mathbb D}$ as the closure of the set consisting of the points $
U: \mathbb R\mapsto \mathbb R^4$ such that there exists one point
$y^i\in \mathbb R$ so that $U-\tilde U\in L^1(\mathbb R, \mathbb
R^4)$ and $TV ( U-\tilde U)\le \delta_0$, where
$$
\tilde U(y)
 =\left\{\begin{array}{ll}U_-,\qquad y<y^i,
 \\U_+,\qquad y^i
\le y\le boundary.
\end{array}\right.
$$
\end{definition}

\begin{remark}
For a solution $U(x, y)$ to the  initial-boundary value problem of
{\rm (\ref{eq1.1})} and {\rm (\ref{IC})}--{\rm (\ref{BC})}, if, for
any fixed $x\ge 0$, $U_x(y)=U(x, y )\in \mathbb D$, then
$y^i=g(0)=0$ when $x=0$, but $y^i<g(x)$ when $x>0$ since there is a
strong shock.
\end{remark}

The semigroup defined by the wave front tracking method is set up in
the following theorem.

\begin{theorem}\label{6.1} If $TV(\bar{U}(\cdot))+ TV(g'(\cdot))$ is
sufficiently small, then  $S^\ep$ defined by the wave front tracking
method is a Cauchy sequence in the $L^1$ sense. Let
$S_x(\bar{U})=\lim_{\ep\rightarrow 0}S_x^\ep(\bar{U})$. There exists
a constant $L$ such that  $S: [0, \infty)\times \mathbb D\mapsto
\mathbb D$ is a uniformly Lipschitz continuous semigroup with
following properties:
\begin{enumerate}
\item[(i)]  $S_0 \bar{U}=\bar{U} $,  $S_{x_1}S_{x_2}\bar{U}=S_{x_1+x_2}\bar{U}$;
\item[(ii)] $||S_x\bar{U}-S_x\bar{V}||_{L^1}\leq L ||\bar{U}-\bar{V}||_{L^1}$;
\item[(iii)] Each trajectory $x\mapsto S_x{\bar U}$ yields an entropy solution
to the
initial-boundary problem {\rm (\ref{eq1.1})} and {\rm
(\ref{IC})}--{\rm (\ref{BC})};
\item[(iv)] If $\bar U\in \mathbb D$ is piecewise constant, then, for $x>0$ sufficiently
small, the function $U(x, \cdot)=S_x{\bar U}$ coincides with the
solution of {\rm (\ref{eq1.1})} and {\rm (\ref{IC})}--{\rm
(\ref{BC})} obtained by piecing together the standard Riemann
solutions and the lateral Riemann solutions.
\end{enumerate}
\end{theorem}

\begin{corollary} If $\, TV(\bar{U}(\cdot))+ TV(g'(\cdot))$ is sufficiently small,
the entropy solution to the initial-boundary problem {\rm
(\ref{eq1.1})} and {\rm (\ref{IC})}--{\rm (\ref{BC})} produced by
the wave front tracking method is unique.
\end{corollary}

To prove Theorem \ref{6.1}, we need the following lemma which can be
found in \cite{BC}.

\begin{lemma}\label{l6.1}
Let $S:[0, \infty)\times \mathbb D\mapsto \mathbb D$ be a globally
Lipschitz semigroup. Let $X>0$,  $\bar {V}\in \mathbb D$, and $V:[0,
X]\mapsto\mathbb D$ be a continuous map whose values are piecewise
constant in the $(x, y)-$plane, with jumps occurring along finitely
many polygonal lines. Let $L$ be the Lipschitz constant of the
semigroup. Then
\begin{equation}
\|V(X)-S_X\bar{V}\|_{L^1}\le L\Big\{\|V(0)-\bar{V}\|_{L^1}
+\int^X_0\overline {\lim_{h\rightarrow 0^+}}
\frac{\|V(x+h)-S_hV(x)\|_{L^1}}{h}dx\Big\}.
\end{equation}
\end{lemma}

\medskip
{\bf Proof of Theorem \ref {6.1}}. It is quite similar to the
proof in \cite{BC}. The difference is that the front tracking method
in \cite{BC} is to use the cut-off function in the order of
$\sqrt{\ep}$, while the front tracking method here is to employ the
simplified Riemann solver rather than the accurate Riemann solver
when the interaction term is less than $\ep$.

By Lemma \ref{l6.1},
\begin{equation}
\|S_X^{\ep_n}\bar{V_n}-S_X^{\ep_m}\bar{V_m}\|_{L^1} \le
L\Big\{\|\bar{V_n}-\bar{V_m}\|_{L^1} +\int^X_0\overline
{\lim_{h\rightarrow 0^+}} \frac{\|S_h^{\ep_n}(S_x^{\ep_m}{\bar
V_m})-S_{x+h}^{\ep_m}{\bar V_m}\|_{L^1}}{h}dx\Big\}.
\end{equation}
Let $\ep_m>\ep_n$. Then the $\ep_m-$approximate solution and
$\ep_n-$approximate solution only differ when the interaction term
of two weak waves or the strength of the wave interacting with the
strong wave is in $[\ep_n,  \ep_m]$. Suppose that there are $N+1$
such interactions. For each weak wave interaction between $\af$ and
$\bt$,
$$
|\alpha\beta|=\ep_m,
$$
and
$$
\mbox{either}\;\alpha\ge \sqrt{\ep_m}\,\,\;\mbox{or}\;\,\, \beta\ge
\sqrt{\ep_m}.
$$
Let $\alpha$ be large. Then
\begin{eqnarray*}\overline{\lim_{h\rightarrow 0^+}}
\frac{\|S_h^{\ep_n}(S_x^{\ep_m}{\bar V_m})
-S_h^{\ep_m}(S_x^{\ep_m}{\bar V_m})\|_{L^1}}{h}
&= & \sum_{i=1}^{N}O(1)\ep_m + O(1)\ep_m\\
&=& \sum_{i=1}^{N}O(1)\sqrt{\ep_m}\alpha + O(1)\ep_m\\
&=& O(1)\sqrt{\ep_m}\,TV(\bar{V}(\cdot))=O(1)\sqrt{\ep_m}.
\end{eqnarray*}
Therefore, $S_x^{\ep_n}\bar {V_n}$ is a Cauchy sequence, which
converges in the $L^1$ sense. Hence, the map $S:[0, \infty)\times
\mathbb D\mapsto \mathbb D$ as the limit of the approximate
solutions produced by the wave front tracking method is
well-defined.

Next, we prove (i) to (iv). Facts (i), (ii), and (iv) are obvious
since $S$ is the limit of $S^\ep$ produced by the wave front
tracking method. It is similar to prove (iii) as \cite{BC}, the only
difference is that the wave front tracking method we employ here is
slightly different. Finally, we can see that the entropy solution
satisfies the boundary condition due to the construction of our
approximate solutions. This completes the proof of Theorem 6.1.

\section{Uniqueness of Entropy solutions in a broader class}

In this section, we first prove the semigroup $S$ defined by the
wave front tracking method is the only standard Riemann semigroup
(SRS) which is defined as Definition \ref{d7.1}. That is, the
semigroup defined by the wave front tracking method is the canonical
trajectory of the standard Riemann semigroup (SRS). Then we prove
that the uniqueness of entropy solutions in a broader class, i.e.
the class of viscosity solutions as defined in \cite{Br2}. The main
point is to prove that, in the class of viscosity solutions, the
entropy solution is unique, which coincides with the trajectory
produced by the wave front tracking method.

\begin{definition}\label{d7.1}
We say that problem {\rm (\ref{eq1.1})} and {\rm (\ref{IC})}--{\rm
(\ref{BC})} admits a standard Riemann semigroup (SRS) if, for some
$\delta_0$, there exists a continuous mapping: $R: [0, \infty)\times
\mathbb D\mapsto \mathbb D$ and a constant $L$ with the following
properties:
\begin{enumerate}
\item[(i)] $R_0 \bar {U}=\bar {U},\; R_{x_1}R_{x_2}\bar{U}=R_{x_1+x_2}\bar{U}$;
\item[(ii)] $\|R_x\bar{U}-R_x\bar{V}\|_{L^1}\le L\|\bar{U}-\bar{V}\|_{L^1}$;
\item[(iii)] If $\bar{U}\in \mathbb D$ is piecewise constant, then,
for $x>0$ sufficiently small, the function $U(x, \cdot)=R_x{\bar U}$
coincides with the solution of {\rm (\ref{eq1.1})} and {\rm
(\ref{IC})}--{\rm (\ref{BC})} obtained by piecing together the
standard Riemann solutions and the lateral Riemann solutions.
\end{enumerate}
\end{definition}

\begin{theorem} \label{t7.1}
Let problem {\rm (\ref{eq1.1})} and {\rm (\ref{IC})}--{\rm
(\ref{BC})} admits a standard Riemann semigroup $R: [0,
\infty)\times \mathbb D\mapsto \mathbb D$. Let $S$ be the semigroup
generated by the wave front tracking method, i.e.
$S_x(\bar{U})=\lim_{\ep\rightarrow 0}S_x^\ep(\bar{U})$.  If $\bar
{U}\in \mathbb D$, then $R_x \bar{U}=S_x \bar{U}$ for all $x\ge 0$.
\end{theorem}

The proof of the theorem is similar to the proof in \cite{Br2} by
using Lemma \ref{l6.1} and the fact that, locally in $x$ direction,
the wave front tracking method and the standard Riemann semigroup
(SRS) both have the structure of the Riemann solutions.

As in \cite{Br2}, there are two types of local approximate
parametrices for {\rm (\ref{eq1.1})}: One is derived from the
self-similar solution of a Riemann problem, and the other is
obtained by ``freezing" the coefficients of the corresponding
quasilinear hyperbolic system in a neighborhood of a given point.

Let $U: [0, \infty)\times \mathbb R\mapsto \mathbb R^4$ be a
function. Fix any point $(\tau, \xi)$ in the domain of $U$. If
$U(\tau, \cdot)\in \mathbb D$, then the bound on the total variation
implies the existence of the limits
$$
U^-=\lim_{y\rightarrow\xi_-}U(\tau, y),
\quad
\;
U^+=\lim_{y\rightarrow\xi_+}U(\tau, y).
$$
Denote by $\omega=\omega (x, y)$ the corresponding solution of the
Riemann problem with $U^-$ and $U^+$ and by $\hat \lambda$  a upper
bound for all characteristic speeds, i.e.,
\begin{equation}\label{lambdahat}
\sup_{U} |\lambda_i(U)|<\hat\lambda, \qquad
i=1,2,3,4.
 \end{equation}

For $x>\tau$,  define the function
\begin{equation}
W^\# _{(U, \tau, \xi)}(x, y)
=\left\{\begin{array}{ll}w(x-\tau, y-\xi)\quad
&\mbox{if}\,|y-\xi|\le\hat \lambda (x-\tau),\\
U(\tau, y) \quad &\mbox{if}\,|y-\xi|>\hat \lambda
(x-\tau).\end{array}\right.
\end{equation}

Set $\tilde A\dot=DW(U(\tau, \xi))$ and $\tilde B\dot=DH(U(\tau,
\xi))$ the Jacobian matrices computed at the point $U(\tau, \xi)$.
For $x>\tau$, define $W^b_{(U, \tau, \xi)}$ as the solution of the
linear Cauchy problem with constant coefficients
 \begin{equation}
 \tilde A V_x + \tilde B V_y=0, \qquad V(\tau, y)=U(\tau, y).
 \end{equation}
Then the functions $W^\#$ and  $W^b$ depend on the values $U(\tau,
\xi)$ and $U(\tau, \xi\pm)$. Next, we define viscosity solutions
which have the same local characterization as $W^\#$ and $W^b$.

\begin {definition} Let $U: [0, X]\mapsto \mathbb D$ be continuous
with respect to the $L^1$ norm. We say  that $U$ is a viscosity
solution of system {\rm (\ref{eq1.1})} and {\rm (\ref{IC})}--{\rm
(\ref{BC})} if there exist constants $C$ and $\hat\lambda$
satisfying {\rm (\ref{lambdahat})} such that, at each point $(\tau,
\xi)\in [0, X)\times \mathbb R$, when $\rho$ and $\ep$ are
sufficiently small,
\begin{eqnarray}
\frac{1}{\ep}\int_{\xi-\rho+\ep\hat\lambda}^{\xi+\rho-\ep\hat\lambda}
&&|U(\tau+\ep, y)-W^\# _{(U, \tau, \xi)}(x, y)|dx\nonumber\\
&\le & C\, TV\{U(\tau)\,:\, (\xi-\rho, \xi)\cup (\xi, \xi+\rho)\},\\
\frac{1}{\ep}\int_{\xi-\rho+\ep\hat\lambda}^{\xi+\rho-\ep\hat\lambda}
&&|U(\tau+\ep, y)-W^b _{(U, \tau, \xi)}(x, y)|dx\nonumber\\
&\le & C\, (TV\{U(\tau)\,:\, (\xi-\rho,  \xi+\rho)\})^2.
\end{eqnarray}
\end{definition}

\begin{lemma} Let $(a, b)$ be a (possibly unbounded) open interval
and let $R$ be the standard Riemann semigroup (SRS). If $\bar {U}$,
$\bar {V}\in\mathbb D$, then, for all $x\ge 0$,
\begin{equation}
\int^{b-\hat\lambda x}_{a+\hat\lambda x}|R_x\bar {U}(y)-R_x\bar
{V}(y)|dy \le L\,\int^b_a |\bar{U}(y)-\bar{V}(y)|dy.\end{equation}
\end{lemma}

\begin{proof} If two initial data functions $\bar {U}$ and $\hat {U}\in\mathbb D$
coincide on $(a, b)$, due to the finite speed of propagation, then
the semigroup generated by the wave front tracking method  has
$S_x\bar{U}(y)=S_x\hat{U}(y)$ for $y\in (a+\hat\lambda x,
b-\hat\lambda x)$.   By Theorem \ref{t7.1},
$R_x\bar{U}(y)=S_x\bar{U}(y)=S_x\hat{U}(y)=R_x\hat{U}(y)$ for all
$y\in (a+\hat\lambda x, b-\hat\lambda x)$. Next, for $\bar {U}$,
$\bar {V}\in\mathbb D$, define
$$
\hat U=\bar{U}\chi_{[a, b]}+ U_-\chi_{(-\infty, a)} + U_+\chi_{(b,
g(x))},
$$
where $\chi$ is the characteristic function, and $U_-$ and $U_+$ are
the states in the definition of $\mathbb D$. It follows that $\hat
U\in \mathbb D$. Similarly, $\hat V=\bar{V}\chi_{[a, b]}+
U_-\chi_{(-\infty, a)} + U_+\chi_{(b, g(x))}\in \mathbb D$. The
uniform Lipschitz continuity of the semigroup $R$ implies
\begin{eqnarray*}
\int^{b-\hat\lambda x}_{a+\hat\lambda x}|R_x\bar {U}(y)-R_x\bar {V}(y)|dy
&=&\int^{b-\hat\lambda x}_{a+\hat\lambda x}|R_x\hat U(y)-R_x\hat V(y)|dy\\
&\le& \|R_x\hat U(y)-R_x\hat V(y)\|\\
&\le& L\, \|\hat U-\hat V\|= L\, \int_a^b|\bar{U}-\bar{V}|dy.
 \end{eqnarray*}
\end{proof}

\begin{theorem}
Assume that problem {\rm (\ref{eq1.1})} and {\rm (\ref{IC})}--{\rm
(\ref{BC})} admits a standard Riemann semigroup $R$. Then a
continuous map $U: [0, X]\mapsto \mathbb D$ is a viscosity solution
of {\rm (\ref{eq1.1})} and {\rm (\ref{IC})}--{\rm (\ref{BC})} if and
only if
\begin{equation}
U(x, \cdot)=R_x \bar {U}\qquad \mbox{for any }x\in [0, X].
\end{equation}
 \end{theorem}

\begin{corollary} For system {\rm (\ref{eq1.1})} and {\rm (\ref{IC})}--{\rm
(\ref{BC})}, the entropy solution is unique in the class of the
viscosity solutions, which coincides with the trajectory
$S_x\bar{U}$ generated by the wave front tracking method, i.e. a
continuous map $U:[0, X]\mapsto \mathbb D$ is a viscosity solution
if and only if
   \begin{equation}
U(x, \cdot)=S_x \bar {U}\qquad \mbox{for any }x\in [0, X].
\end{equation}
 \end{corollary}

The proof is similar to the argument in \cite{Br2}. The only
difference is that there is a strong shock in our case; however, we
can still carry out the proof as long as the convergence of the wave
front tracking method is achieved which has been proved in Section
3.
\begin{remark}
For the potential flow, isentropic or isothermal Euler flow {\rm
(\ref{eq1.2})}, which are the simpler cases as the $L^1$ stability
problem as concerned, we obtain the same results as the full Euler
equations.
\end{remark}

\bigskip
\noindent{\bf Acknowledgments.} The research of Gui-Qiang Chen was
supported in part by the National Science Foundation under Grants
DMS-0505473, DMS-0426172, DMS-0244473, and an Alexandre von Humboldt
Foundation Fellowship. The research of Tian-Hong Li was supported in
part by the National Science Foundation under Grants DMS-0244383 and
DMS-0244473. The authors would like to thank Professor Tai-Ping Liu
for helpful discussion.

\end{document}